\documentclass[11pt]{amsart}
\usepackage{amsmath}
\usepackage{amssymb}
\usepackage{gastex}
\usepackage{ifthen}
\usepackage{url}

\newcommand{\mediumSize}[1]{\fontsize{9pt}{12pt}\selectfont #1\normalsize}
\newcommand{\mediumFont}[1]{\normalfont\mediumSize{#1}}
\newcommand{\malcev}%
  {\mathop{\text{\normalsize{\raisebox{0.1mm}{\textcircled{\raisebox{0.2mm}{\mediumFont{\textit{m}}}}}}}}}
\font\petite=cmmi10 at 8pt
\def\malcev{\mathbin{\hbox{$\bigcirc$\rlap{\kern-9pt\raise0,75pt\hbox{\petite m}}}}}

\newcommand{\varietyFont}[1]{\mathrm{\mathbf{#1}}}
\newcommand{\DA}{\varietyFont{D\hspace{-1pt}A}}
\newcommand{\LI}{L\varietyFont{I}}

\newcommand{\Jone}{\varietyFont{SL}}

\newcommand{\DS}{\varietyFont{DS}}

\let\J\Jtrivial
\newcommand{\RR}{\varietyFont{R}}
\newcommand{\LL}{\varietyFont{L}}

\newcommand{\I}{\varietyFont{I}}
\newcommand{\Nil}{\varietyFont{Nil}}
\newcommand{\K}{\varietyFont{K}}
\newcommand{\D}{\varietyFont{D}}
\newcommand{\B}{\varietyFont{B}}
\newcommand{\BL}{\varietyFont{BL}}
\newcommand{\BR}{\varietyFont{BR}}
\newcommand{\Y}{\varietyFont{Y}}
\newcommand{\RZ}{\varietyFont{RZ}}
\newcommand{\LZ}{\varietyFont{LZ}}

\newcommand{\M}{\varietyFont{M}}
\newcommand{\V}{\varietyFont{V}}  
\newcommand{\W}{\varietyFont{W}}
\newcommand{\X}{\varietyFont{X}}

\newcommand{\LangV}{\mathcal{V}}  

\def\calB{\mathcal{B}}
\def\calJ{\mathcal{J}}
\def\calL{\mathcal{L}}

\newcommand{\greenFont}[1]{\mathcal{#1}}
\newcommand{\gR}{\mathrel{\greenFont{R}}}
\newcommand{\gL}{\mathrel{\greenFont{L}}}
\newcommand{\gJ}{\mathrel{\greenFont{J}}}

\theoremstyle{plain}

\newtheorem{theorem}{Theorem}[section]
\newtheorem{proposition}[theorem]{Proposition}

\newtheorem{corollary}[theorem]{Corollary}
\newenvironment{example}[1][]{\ifthenelse{\equal{#1}{}}{\begin{expl}\upshape}{\begin{expl}[#1]\upshape}}{\hspace*{\fill}$\Diamond$\end{expl}}

\theoremstyle{definition}

\newtheorem{remark}[theorem]{Remark}
\newtheorem{expl}[theorem]{Example}

\newenvironment{Proof}{\rm \trivlist \item[\hskip \labelsep{\bf
Proof.}]}{\cqfd\endtrivlist}

\def\cqfd{\skip10=\parfillskip\parfillskip=0pt
\enspace\hfill\symbolecqfd\par\parfillskip=\skip10\par\medskip}
\def\symbolecqfd{\rlap{$\sqcap$}$\sqcup$}

\def\proof{\begin{Proof}}
\def\eop{\end{Proof}}

\newcounter{commentcounter}

\def\inv{^{-1}}
\let\phi\varphi
\let\epsilon\varepsilon

\def\llbracket{[\![}
\def\rrbracket{]\!]}

\begin{document}

\title{On the lattice of sub-pseudovarieties of $\DA$}
\author{Manfred Kufleitner
\and
Pascal Weil}
\thanks{{}\\
\url{manfred.kufleitner@fmi.uni-stuttgart.de}
\ Universit\"at Stuttgart\\
Institut f\"ur Formale Methoden, Universit\"atsstr.  38, D-70569
Stuttgart, Germany.\\
\url{pascal.weil@labri.fr}\ Universit\'e de Bordeaux, CNRS and IIT Delhi\\
LaBRI, 351 cours de la Lib\'eration, 33405 Talence Cedex, France.}
%
%
\date{6 April 2009}
\begin{abstract}
    The wealth of information that is available on the lattice of
    varieties of bands, is used to illuminate the structure of the
    lattice of sub-pseudovarieties of $\DA$, a natural generalization
    of bands which plays an important role in language theory and in
    logic.  The main result describes a hierarchy of decidable
    sub-pseudovarieties of $\DA$ in terms of iterated Mal'cev products
    with the pseudovarieties of definite and reverse definite semigroups.
\end{abstract}

\maketitle


The complete elucidation of the structure of the lattice $\calL\calB$
of band varieties is one of the jewels of semigroup theory: this
lattice turns out to be countable, with a simple structure (Birjukov
\cite{Birjukov1970al}, Fennemore \cite{Fennemore1971mn}, Gerhard
\cite{Gerhard1970ja}, see Section~\ref{lattice of bands} below for the
main features of this structure).  Moreover, each of its elements can
be defined by a small number of identities (at most 3), and we can
efficiently solve the membership problem in each variety of bands, as
well as the word problem in its free object
\cite{GerhardPetrich1989plms}.

As bands are locally finite, the lattice $\calL(\B)$ of
pseudovarieties of finite bands is isomorphic to $\calL\calB$: a class
of finite bands is a pseudovariety if and only if it is the class of
finite elements of a variety.

In this paper, we discuss the structure of the lattice of
sub-pseudovarieties of $\DA$, which is a natural generalization of the
pseudovariety $\B$ of bands.  Indeed, $\DA$ is the maximum
pseudovariety in which all regular elements are bands.  This
pseudovariety actually has several other interesting algebraic
characterizations, and also many other characterizations in terms of
formal languages and logic, see the survey by Tesson and Th\'erien
\cite{TessonTherien2002}.  This only adds to the motivation to better
understand the lattice $\calL(\DA)$ of its sub-pseudovarieties.

In fact the authors' initial motivation regards one of the logical
characterizations of $\DA$ by means of the 2-variable fragment of
first-order theory of the linear order \cite{TherienWilke1998stoc},
and the main result of this paper finds an application in a paper on
the language-theoretic characterizations of the quantifier alternation
hierarchy within that logic \cite{KufleitnerWeilFO2mfcs}.  The
characterization in \cite{KufleitnerWeilFO2mfcs} can be viewed as an
algebraic counterpart of recent results of Weis and Immerman's
description on the 2-variable fragment of first-order logic
\cite{WeisImmerman2009lmcs}, and the main result of the present paper 
gives a purely algebraic foundation to the results in 
\cite{WeisImmerman2009lmcs} and \cite{KufleitnerWeilFO2mfcs}.

Trotter and Weil \cite{TrotterWeil1997au} initiated the study of the
structure of $\calL(\DA)$ by considering the map $\V \mapsto
\V\cap\B$, from $\calL(\DA)$ to $\calL(\B)$.  They showed that, for
each pseudovariety of bands $\Y$, the inverse image of $\Y$ is an
interval in $\calL(\DA)$, with minimum element $\Y$ itself.  They also
showed how to effectively turn the identities defining $\Y$ as a band
pseudovariety, into pseudo-identities defining $\Y^\uparrow$, the
maximal element of that interval.  This result uncovers the
interesting role played by the lattice of decidable pseudovarieties
given by the $\Y^\uparrow$, $\Y\in\calL(\B)$.

The missing element was an understanding of the fashion in which one
can \textit{climb} in that lattice.  The beautiful results on
$\calL(\B)$ include a description of the different levels of the
hierarchy it forms, in terms of Mal'cev products with the
pseudovarieties $\RZ$ and $\LZ$ of right zero and left zero bands.  In
this paper, we elucidate the structure of the sublattice of
$\calL(\DA)$ formed by the $\Y^\uparrow$ ($\Y \in \calL(\B)$), in
terms of Mal'cev products as well, with definite and reverse definite
semigroups.  This helps establish that the $\Y^\uparrow$ form an
infinite hierarchy, whose union is all of $\DA$.  It follows in
particular that $\DA$ is the least pseudovariety containing
semilattices, which is closed under Mal'cev product with definite and
reverse definite semigroups, -- a fact with an interesting interpretation in
formal language theory.

Interestingly, this last result was recently proved, independently 
and by completely different means (logical and language theoretical) 
by Lodaya, Pandya and Shah \cite{LodayaPS2008ifip}.

The paper is organized as follows: Section~\ref{sec: prelim psv}
summarizes what the reader needs to know (for the purpose of this
paper!)  about pseudovarieties and Mal'cev products.
Section~\ref{sec: prelim bands and DA} discusses the known results on
bands, $\DA$ and their respective lattice of sub-pseudovarieties, and
Section~\ref{sec: main result} gives our main result.  Its
consequences are discussed in Section~\ref{sec: applications}, in
semigroup- and in language-theoretic terms.

\section{Preliminaries on pseudovarieties}\label{sec: prelim psv}

\subsection{Pseudovarieties}

Recall that a pseudovariety of semigroups (resp.  monoids) is a class
of finite semigroups (resp.  monoids) closed under taking quotients,
finite direct products and subsemigroups (resp.  submonoids). If $\V$ 
is a pseudovariety of semigroups, we denote by $\V_\M$ the 
pseudovariety of monoids which consists of the monoids in $\V$. The 
pseudovariety $\V$ is called \textit{monoidal} if it is generated by 
the monoids it contains.

If $\W$ is a pseudovariety of monoids, we denote by $L\W$ the class of
semigroups $S$ such that, for each idempotent $e$, the monoid $eSe \in
\W$: $L\W$ forms a pseudovariety, the largest one such that $(L\W)_\M
\subseteq \W$.

There is a vast literature on pseudovarieties, and on their definition
by \textit{pseudo-identities}, see
\cite{Almeida1994book,Weil2002ijac}.  For our purpose, it is enough to
consider so-called \textit{$\omega$-pseudo-identities} of the form $u
= v$, where $u$ and $v$ are obtained from a countable alphabet of
symbols $X$ using the operation of concatenation and formal
$(\omega-1)$-power.  For instance, we will consider in the sequel
identities like $(xy)^\omega x (xy)^\omega = (xy)^\omega$ -- where
$z^\omega$ stands for $z^{\omega-1}z$.  A finite semigroup $S$
satisfies the $\omega$-pseudo-identity $u=v$ if, for every map
$\phi\colon X \to S$, we have $\hat\phi(u) = \hat\phi(v)$, where
$\hat\phi$ extends $\phi$ to a monoid morphism such that
$\hat\phi(t^\omega)$ is the (unique) idempotent power of $\hat\phi(t)$
and $\hat\phi(t^{\omega-1})$ is the inverse of the element
$\hat\phi(t^\omega)\hat\phi(t)$ in the minimal ideal of the
subsemigroup generated by $\hat\phi(t)$, which is a group
\cite{Almeida1994book,Pin1986book}.  By a common abuse of notation, we
also denote by $s^\omega$ ($s\in S$) the idempotent power of $s$, and
by $s^{\omega-1}$ the inverse of $s^\omega s$ in the minimal ideal of
the subsemigroup generated by $s$.

If $(u_i = v_i)_{i\in I}$ is a family of $\omega$-pseudo-identities,
we denote by $\llbracket (u_i = v_i)_{i\in I} \rrbracket$ the class of
finite semigroups which satisfy each $\omega$-pseudo-identity $u_i =
v_i$.  Such a class is always a pseudovariety\footnote{Not all
pseudovarieties are obtained this way; for a more rigorous discussion
of pseudidentities, and in particular for a converse statement
(involving a much larger set of pseudo-identities), see
\cite{Almeida1994book,Weil2002ijac}.}.

\subsection{Mal'cev products}\label{sec: Malcev}

Let $\V$ be a pseudovariety of semigroups, $\W$ a pseudovariety of
semigroups (resp.  monoids) and $M$ a finite semigroup (resp.
monoid).  We say that $M \in \V\malcev\W$ (the \textit{Mal'cev
product} of $\V$ and $\W$) if there exists a finite monoid $T$ and
onto morphisms $\alpha\colon T\rightarrow M$ and $\beta\colon
T\rightarrow N$ such that $N\in \W$ and, for each idempotent $e$ of
$N$, $\beta\inv(e) \in \V$ (we say that $\beta$ is a
\textit{$\V$-morphism)}.  Then $\V\malcev\W$ is a pseudovariety of
semigroups (resp.  monoids), see
\cite{Pin1986book,Almeida1994book,PinWeil1996ja}.

In the sequel, we will consider Mal'cev products where the first 
component is one of the pseudovarieties $\Nil$, $\LZ$, $\RZ$, $\K$, $\D$ and 
$L\I$, which are defined as follows:
\begin{align*}
    &\K = \llbracket x^\omega y = x^\omega\rrbracket, \qquad \D = 
    \llbracket yx^\omega = x^\omega\rrbracket, \\
    &\Nil = \K \cap \D = \llbracket x^\omega y = yx^\omega =
    z^\omega\rrbracket, \\
    &L\I = \K \vee \D = \llbracket x^\omega yx^\omega =
    x^\omega\rrbracket, \\
    &\LZ = \K\cap \llbracket x^2 = x\rrbracket = \llbracket xy = x,
    x^2 = x\rrbracket,\\
    &\RZ = \D\cap\llbracket x^2 = x\rrbracket = \llbracket yx = x, x^2
    = x\rrbracket \\
\end{align*}
We will use the following fact, due to Krohn, Rhodes and Tilson
\cite{KrohnRT1965arbib}, see \cite[Corollary 4.3]{HallWeil1999sf}.
The $\gJ$-quasi-order on $M$ is defined as follows: $x \le_{\gJ} y$ if
and only if $x = uyv$ for some $u,v \in M\cup\{1\}$. We write $x 
<_{\gJ} y$ if $x \le_{\gJ} y$ but not $y \le_{\gJ} x$; that is, if 
the 2-sided ideal of $M$ generated by $x$ is properly contained in 
the 2-sided ideal generated by $y$.

\begin{proposition}\label{congruence malcev}
    Let $M$ be a finite semigroup and let $\sim_\K$ and $\sim_\D$ be
    the equivalence relations $\sim_\K$ and $\sim_\D$ on $M$ given,
    for $s,t \in M$, by
    \begin{align*}
	s\sim_\K t &\hbox{ if and only if, for all } e\in E(M),\ es,et
	<_{\gJ} e \hbox{ or } es = et \\
	s\sim_\D t &\hbox{ if and only if, for all } e\in E(M),\ se,te
	<_{\gJ} e \hbox{ or } es = et.
    \end{align*}
    These two relations are congruences and $M/\!\sim_\K$ 
    (resp. $M/\!\sim_\D$) is the least quotient of $M$ such that the 
    projection is a $\K$- (resp. a $\D$-) morphism.
    
    If $\V$ is a pseudovariety of semigroups (resp.  monoids) and $M$
    is a finite semigroup (resp.  monoid), then $M \in \K\malcev\V$
    (resp.  $M \in \D\malcev\V$) if and only if $M/\!\sim_\K\in \V$
    (resp.  $M/\!\sim_\D \in \V$).

    In particular, if $\V$ is decidable, then so are $\K\malcev\V$ and
    $\D\malcev\V$.
\end{proposition}

\section{Preliminaries on bands and $\DA$}\label{sec: prelim bands 
and DA}

\subsection{Bands and $\DA$}

A \textit{band} is a semigroup in which every element is idempotent.
We denote by $\B$ the pseudovariety of bands, that is, $\B =
\llbracket x^2 = x\rrbracket = \llbracket x^\omega = x\rrbracket$.

Let $\DA = \llbracket (xy)^\omega x (xy)^\omega = (xy)^\omega
\rrbracket$.  The following result combines several known results: we
refer the reader to \cite{TessonTherien2002} for a synthesis on $\DA$ 
(see also  
\cite{Schutzenberger1976sf,FichBrzozowksi1979icalp,Weil1993sf}).

\begin{proposition}\label{prop basics DA}
    If $M$ is a finite semigroup, the following are equivalent.
    \begin{itemize}
	\item[(1)] $M \in \DA$,
	
	\item[(2)] every regular element of $M$ is idempotent,
	
	\item[(3)] for every idempotent $e\in E(M)$, we have
	$eM_ee = e$, where $M_e$ is the sub-semigroup of $M$
	generated the elements $x\ge_\calJ e$,
	
	\item[(4)] $M \in \LI\malcev \Jone$, where $\Jone = \llbracket
	x^2 = x, xy = yx\rrbracket$ is the pseudovariety of idempotent
	and commutative semigroups.
    \end{itemize}
\end{proposition}

\begin{corollary}\label{DA max}
    $\DA$ (resp.  $\DA_\M$) is the maximum pseudovariety of
    semigroups (resp. monoids), in which every regular semigroup 
    (resp. monoid) is a band.
\end{corollary}

\proof
If $M\in \DA$ and a regular semigroup, then $M$ is a band by 
Proposition~\ref{prop basics DA} (2).

Let now $\V$ be a pseudovariety of semigroups in which every regular
element is a band.  Recall that an element $s\in M$ is
\textit{regular} if there exists $t\in M$ such that $sts = s$.  In
particular, the regular elements of $M$ are exactly the elements of
the form $(st)^\omega s$ ($s,t \in M$). Therefore $\V$ satisfies the 
$\omega$-pseudo-identity $(xy)^\omega x (xy)^\omega x = (xy)^\omega 
x$, and by right multiplication by $y(xy)^{\omega-1}$, we find that 
$\V$ satisfies $(xy)^\omega x (xy)^\omega = (xy)^\omega$, that is, 
$\V$ is contained in $\DA$.
\eop

Many more characterizations of $\DA$ can be found in the literature,
see \cite{TessonTherien2002,TessonTherien2007lmcs}.  In this paper, we
will encounter one more, in Section~\ref{sec: languages} below, in
relation with formal language theory.

\subsection{The lattice $\calL(\B)$}\label{lattice of bands}

Since the free band over a finite alphabet is finite (see
\cite{Howie1995book}), the lattice $\calL(\B)$ of sub-pseudovarieties
of $\B$ is isomorphic with the lattice $\calL\calB$ of all band
varieties.  The structure of that lattice was elucidated around 1970
(Birjukov \cite{Birjukov1970al}, Fennemore
\cite{Fennemore1971mn}, Gerhard \cite{Gerhard1970ja}).  The lattice
$\calL(\B)$ turns out to be countable, with a simple structure.  We
summarize below the main results concerning this lattice that will be
useful to us.

We define the pseudovarieties $\BR_m$, $\BR'_{m+1}$, $\BL_m$ and
$\BL'_{m+1}$ ($m \ge 1$) by letting\footnote{In the traditional
terminology of bands, the elements of $\BR'_2$ and $\BL'_2$ are called
\textit{right normal} and \textit{left normal} bands respectively.}
\begin{align*}
    &\BR_1 = \BL_1 = \Jone,&\\
    &\BR_{m+1} = \LZ\malcev\BL_m,&\BL_{m+1} = \RZ\malcev\BR_m,\\
    &\BR'_2 = \llbracket xyz = xzy, x^2 = x\rrbracket,&\BL'_2 =
    \llbracket xyz = yxz, x^2 = x\rrbracket,\\
    &\BR'_{m+2} = \LZ\malcev\BL'_{m+1},&\BL'_{m+2} =
    \RZ\malcev\BR'_{m+1}.
\end{align*}
The following statement describes the structure of the lattice
$\calL(\B)$, as discussed by Gerhard and Petrich \cite{GerhardPetrich1989plms}
(the last item is due to Wismath \cite{Wismath1986sf}).

\begin{figure}
    \begin{picture}(63,130)(0,-130)

    \node[Nfill=y,ExtNL=y,NLangle=0.0,NLdist=2.0,Nw=1.0,Nh=1.0,Nmr=0.5](I)(32.0,-130.0){$\I$}
    \node[ExtNL=y,NLangle=180.0,NLdist=2.0,Nw=1.0,Nh=1.0,Nmr=0.5](LZ)(16.0,-118.0){$\LZ$}
    \node[ExtNL=y,NLangle=0.0,NLdist=2.0,Nw=1.0,Nh=1.0,Nmr=0.5](RZ)(48.0,-118.0){$\RZ$}
    \node[Nw=1.0,Nh=1.0,Nmr=0.5](ReB)(32.0,-106.0){}

    \node[Nfill=y,ExtNL=y,NLangle=0.0,NLdist=2.0,Nw=1.0,Nh=1.0,Nmr=0.5](J1)(32.0,-118.0){$\Jone$}
    \node[ExtNL=y,NLangle=180.0,NLdist=2.0,Nw=1.0,Nh=1.0,Nmr=0.5](R'0)(16.0,-106.0){$\BR'_2$}
    \node[ExtNL=y,NLangle=0.0,NLdist=2.0,Nw=1.0,Nh=1.0,Nmr=0.5](L'0)(48.0,-106.0){$\BL'_2$}

    \node[Nfill=y,ExtNL=y,NLangle=180.0,NLdist=2.0,Nw=1.0,Nh=1.0,Nmr=0.5](R1)(0.0,-94.0){$\BR_2$}
    \node[Nfill=y,ExtNL=y,NLangle=0.0,NLdist=2.0,Nw=1.0,Nh=1.0,Nmr=0.5](L1)(64.0,-94.0){$\BL_2$}
    \node[ExtNL=y,NLangle=0.0,NLdist=2.0,Nw=1.0,Nh=1.0,Nmr=0.5](M1)(32.0,-94.0){}

    \node[ExtNL=y,NLangle=180.0,NLdist=2.0,Nw=1.0,Nh=1.0,Nmr=0.5](N1)(16.0,-82.0){}
    \node[ExtNL=y,NLangle=0.0,NLdist=2.0,Nw=1.0,Nh=1.0,Nmr=0.5](P1)(48.0,-82.0){}

    \node[ExtNL=y,NLangle=180.0,NLdist=2.0,Nw=1.0,Nh=1.0,Nmr=0.5](R'1)(0.0,-70.0){$\BR'_3$}
    \node[ExtNL=y,NLangle=0.0,NLdist=2.0,Nw=1.0,Nh=1.0,Nmr=0.5](L'1)(64.0,-70.0){$\BL'_3$}
    \node[Nfill=y,ExtNL=y,NLangle=0.0,NLdist=2.0,Nw=1.0,Nh=1.0,Nmr=0.5](M'1)(32.0,-70.0){}

    \node[ExtNL=y,NLangle=180.0,NLdist=2.0,Nw=1.0,Nh=1.0,Nmr=0.5](N'1)(16.0,-58.0){}
    \node[ExtNL=y,NLangle=0.0,NLdist=2.0,Nw=1.0,Nh=1.0,Nmr=0.5](P'1)(48.0,-58.0){}

    \node[Nfill=y,ExtNL=y,NLangle=180.0,NLdist=2.0,Nw=1.0,Nh=1.0,Nmr=0.5](R2)(0.0,-46.0){$\BR_3$}
    \node[Nfill=y,ExtNL=y,NLangle=0.0,NLdist=2.0,Nw=1.0,Nh=1.0,Nmr=0.5](L2)(64.0,-46.0){$\BL_3$}
    \node[ExtNL=y,NLangle=0.0,NLdist=2.0,Nw=1.0,Nh=1.0,Nmr=0.5](M2)(32.0,-46.0){}

    \node[ExtNL=y,NLangle=180.0,NLdist=2.0,Nw=1.0,Nh=1.0,Nmr=0.5](N2)(16.0,-34.0){}
    \node[ExtNL=y,NLangle=0.0,NLdist=2.0,Nw=1.0,Nh=1.0,Nmr=0.5](P2)(48.0,-34.0){}

    \node[ExtNL=y,NLangle=180.0,NLdist=2.0,Nw=1.0,Nh=1.0,Nmr=0.5](R'2)(0.0,-22.0){$\BR'_4$}
    \node[ExtNL=y,NLangle=0.0,NLdist=2.0,Nw=1.0,Nh=1.0,Nmr=0.5](L'2)(64.0,-22.0){$\BL'_4$}
    \node[Nfill=y,ExtNL=y,NLangle=0.0,NLdist=2.0,Nw=1.0,Nh=1.0,Nmr=0.5](M'2)(32.0,-22.0){}

    \node[ExtNL=y,NLangle=180.0,NLdist=2.0,Nw=1.0,Nh=1.0,Nmr=0.5](N'2)(16.0,-10.0){}
    \node[ExtNL=y,NLangle=0.0,NLdist=2.0,Nw=1.0,Nh=1.0,Nmr=0.5](P'2)(48.0,-10.0){}

    \node[Nfill=y,ExtNL=y,NLangle=0.0,NLdist=2.0,Nw=0.0,Nh=0.0,Nmr=0.5](topR)(56.0,-4.0){}
    \node[Nfill=y,ExtNL=y,NLangle=0.0,NLdist=2.0,Nw=0.0,Nh=0.0,Nmr=0.5](topL)(8.0,-4.0){}
    \node[Nfill=y,ExtNL=y,NLangle=0.0,NLdist=2.0,Nw=0.0,Nh=0.0,Nmr=0.5](middleL)(40.0,-4.0){}
    \node[Nfill=y,ExtNL=y,NLangle=0.0,NLdist=2.0,Nw=0.0,Nh=0.0,Nmr=0.5](middleR)(24.0,-4.0){}

    \drawedge[dash={1.0 1.0 1.0 1.0}{0.0},AHnb=0](I,LZ){ }
    \drawedge[dash={1.0 1.0 1.0 1.0}{0.0},AHnb=0](I,RZ){ }
    \drawedge[dash={1.0 1.0 1.0 1.0}{0.0},AHnb=0](RZ,L'0){ }
    \drawedge[dash={1.0 1.0 1.0 1.0}{0.0},AHnb=0](LZ,R'0){ }
    \drawedge[dash={1.0 1.0 1.0 1.0}{0.0},AHnb=0](RZ,ReB){ }
    \drawedge[dash={1.0 1.0 1.0 1.0}{0.0},AHnb=0](LZ,ReB){ }
    \drawedge[dash={1.0 1.0 1.0 1.0}{0.0},AHnb=0](ReB,M1){ }

    \drawedge[dash={1.0 1.0 1.0 1.0}{0.0},AHnb=0](R'0,M1){ }
    \drawedge[dash={1.0 1.0 1.0 1.0}{0.0},AHnb=0](M1,P1){ }
    \drawedge[dash={1.0 1.0 1.0 1.0}{0.0},AHnb=0](P1,L'1){ }
    \drawedge[dash={1.0 1.0 1.0 1.0}{0.0},AHnb=0](L'0,M1){ }
    \drawedge[dash={1.0 1.0 1.0 1.0}{0.0},AHnb=0](M1,N1){ }
    \drawedge[dash={1.0 1.0 1.0 1.0}{0.0},AHnb=0](N1,R'1){ }
    \drawedge[dash={1.0 1.0 1.0 1.0}{0.0},AHnb=0](R'1,N'1){ }
    \drawedge[dash={1.0 1.0 1.0 1.0}{0.0},AHnb=0](N'1,M2){ }
    \drawedge[dash={1.0 1.0 1.0 1.0}{0.0},AHnb=0](M2,P2){ }
    \drawedge[dash={1.0 1.0 1.0 1.0}{0.0},AHnb=0](P2,L'2){ }
    \drawedge[dash={1.0 1.0 1.0 1.0}{0.0},AHnb=0](L'1,P'1){ }
    \drawedge[dash={1.0 1.0 1.0 1.0}{0.0},AHnb=0](P'1,M2){ }
    \drawedge[dash={1.0 1.0 1.0 1.0}{0.0},AHnb=0](M2,N2){ }
    \drawedge[dash={1.0 1.0 1.0 1.0}{0.0},AHnb=0](N2,R'2){ }
    \drawedge[dash={1.0 1.0 1.0 1.0}{0.0},AHnb=0](R'2,N'2){ }
    \drawedge[dash={1.0 1.0 1.0 1.0}{0.0},AHnb=0](N'2,middleR){ }
    \drawedge[dash={1.0 1.0 1.0 1.0}{0.0},AHnb=0](L'2,P'2){ }
    \drawedge[dash={1.0 1.0 1.0 1.0}{0.0},AHnb=0](P'2,middleL){ }

    \drawedge[AHnb=0](I,J1){ }
    \drawedge[AHnb=0](J1,R'0){ }
    \drawedge[AHnb=0](R'0,R1){ }
    \drawedge[AHnb=0](J1,L'0){ }
    \drawedge[AHnb=0](L'0,L1){ }
    \drawedge[AHnb=0](L1,P1){ }
    \drawedge[AHnb=0](P1,N'1){ }
    \drawedge[AHnb=0](N'1,R2){ }
    \drawedge[AHnb=0](R1,N1){ }
    \drawedge[AHnb=0](N1,P'1){ }
    \drawedge[AHnb=0](P'1,L2){ }
    \drawedge[AHnb=0](R2,N2){ }
    \drawedge[AHnb=0](N2,P'2){ }
    \drawedge[AHnb=0](P'2,topR){ }
    \drawedge[AHnb=0](L2,P2){ }
    \drawedge[AHnb=0](P2,N'2){ }
    \drawedge[AHnb=0](N'2,topL){ }

    \end{picture}

    \caption{The lattice $\calL(\B)$; solid lines and bullets denote
    the monoidal pseudovarieties}\label{figure lattice}
\end{figure}
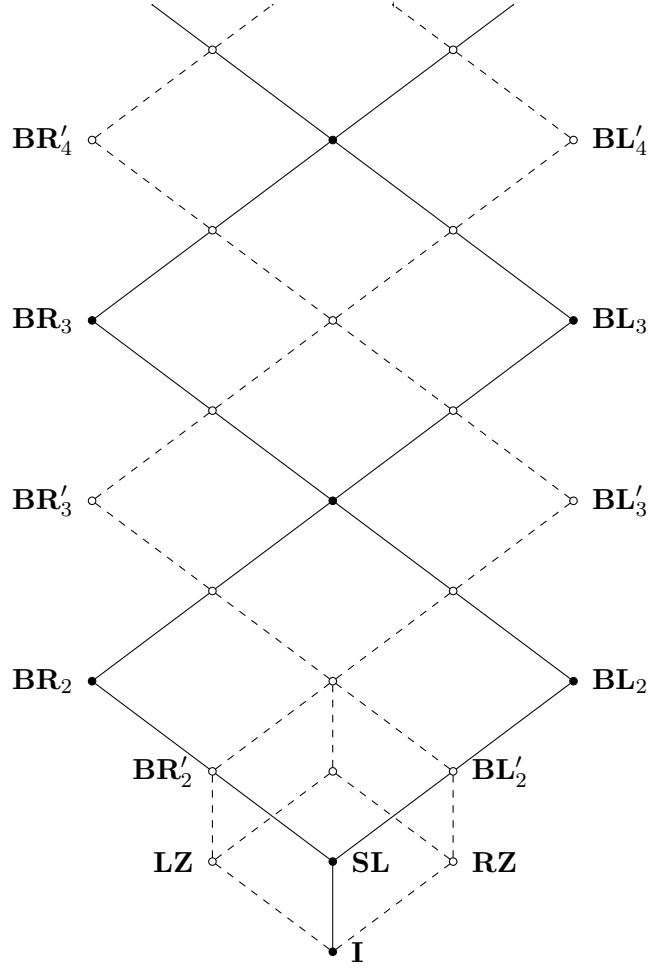

\begin{theorem}\label{lattice bands}
    \begin{itemize}
	\item[(1)] The lattice $\calL(\B)$ consists of the trivial
	pseudovariety $\I$, the pseudovariety $\B$, the $\BR_m$,
	$\BR'_{m+1}$, $\BL_m$, $\BL'_{m+1}$ ($m \ge 1$) and their
	intersections.  It is depicted in Figure~\ref{figure lattice}
	(omitting its top element $\B$).
	
	\item[(2)] The monoidal band pseudovarieties are the trivial
	pseudovariety $\I$, $\B$, and the $\BR_m$, $\BL_m$ and $\BR_m
	\cap \BL_m$ ($m \ge 1$).
	
	\item[(3)] For each $m\ge 1$, we have $\BR_m \vee \BL_m =
	\BR_{m+1}\cap \BL_{m+1}$.
	
	\item[(4)] For each $m\ge 2$, $\BR'_{m+1} = L\BR_m \cap \B$
	and $\BL'_{m+1} = L\BL_m \cap \B$.
	
	\item[(5)] For each $m\ge 2$, the $m$-generated free band lies
	in $\BR'_m \vee \BL'_m = \BR'_{m+1} \cap \BL'_{m+1}$.
	
	\item[(6)] $\calL(\B_\M)$ is isomorphic to the lattice of
	monoidal band pseudovarieties, with the isomorphism given by
	$\V \mapsto \V_\M$ ($\V\in \calL(\B)$, monoidal).
    \end{itemize}
\end{theorem}

Gerhard and Petrich \cite{GerhardPetrich1989plms} also give identities
defining the band pseudovarieties.  Let $x_1,x_2,\dots$ be a sequence
of variables.  If $u$ is a word on that alphabet, we let $\bar u$ be
the mirror image of $u$, that is, the word obtained from reading $u$
from right to left.  We let
\begin{align*}
    &G_2 = x_2x_1,\qquad I_2 = x_2x_1x_2\\
    \hbox{and for $m \ge 2$\quad}&G_{m+1} = x_{m+1}\overline{G_m},\quad I_{m+1} =
    G_{m+1}x_{m+1}\overline{I_m}.
\end{align*}

\begin{theorem}\label{equations bands}
    For each $m\ge 2$, we have $\BR_m = \llbracket x^2 = x, G_m =
    I_m\rrbracket$ and $\BL_m = \llbracket x^2 = x, \overline{G_m} =
    \overline{I_m}\rrbracket$.
\end{theorem}

Note that Theorems~\ref{lattice bands} and~\ref{equations bands} allow
the computation of defining identities for each band pseudovariety,
and indeed to show that each can be defined by a set of at most three
identities.

\subsection{The map $\V \mapsto \V\cap\B$}\label{sec: projecting lattices}

The map $\V \mapsto \V\cap\B$ from $\calL(\DA)$ to $\calL(\B)$ can be
used to derive information on $\calL(\DA)$ from the information
available on $\calL(\B)$.  The following statement was proved by
Trotter and Weil \cite{TrotterWeil1997au} (and, independently, by
Reilly and Zhang \cite{ReillyZhang1997bams} for the first item).

\begin{theorem}\label{thm: TW}
    \begin{itemize}
	\item[(1)] The map $\V \mapsto \V\cap\B$ from $\calL(\DA)$ to
	$\calL(\B)$ is a complete lattice morphism, and the inverse
	image of a band pseudovariety $\Y$ is an interval of
	the form $[\Y,\Y^\uparrow]$.
	
	\item[(2)] For each $\Y\in \calL(\B)$, we have $\Y^\uparrow = 
	(\LZ \malcev \Y)^\uparrow \cap (\RZ \malcev \Y)^\uparrow$.
	
	\item[(3)] The mapping $\V \mapsto \V\cap\B_\M$ from
	$\calL(\DA_\M)$ to $\calL(\B_\M)$ shares the properties from
	statement (1).  Moreover, if $\Y$ is a monoidal band
	pseudovariety, then $(\Y_\M)^\uparrow = (\Y^\uparrow)_\M$.
    
	\item[(4)] If $\Y$ is a monoidal pseudovariety of bands, then
	$(L\Y_\M \cap \B)^\uparrow = L(\Y^\uparrow_\M) \cap \DA$.
	
	\item[(5)] For each $m\ge 2$, $\BR_m^\uparrow = \DA \cap
	\llbracket \phi(G_m) = \phi(I_m)\rrbracket$ and
	$\BL_m^\uparrow = \DA \cap \llbracket \phi(\overline{G_m}) =
	\phi(\overline{I_m})\rrbracket$, where $\phi$ is given by
	\begin{align*}
	    &\phi(x_1) = (x_1^\omega x_2^\omega
	    x_1^\omega)^\omega,\quad \phi(x_2) = x_2^\omega\\
	    \hbox{ and, for $m \ge 2$, }&\phi(x_{m+1}) = (x_{m+1}^\omega
	    \phi(\overline{G_m}G_m)^\omega x_{m+1}^\omega)^\omega.
	\end{align*}
	
	\item[(6)] Let $m\ge 1$.  Every $m$-generated semigroup in
	$\DA$ is in $\Jone^\uparrow$ if $m = 1$, in
	$\BR_3^\uparrow\cap\BL_3^\uparrow$ if $m=2$, and in
	$L(\BR_m^\uparrow\cap\BL_m^\uparrow)$ if $m\ge 3$.  In every
	case, such a semigroup is in $\BR_{m+1}^\uparrow \cap
	\BL_{m+1}^\uparrow$
    \end{itemize}
\end{theorem}

It is elementary to verify \cite{TrotterWeil1997au} that $\I^\uparrow
= \Nil$, $\LZ^\uparrow = \K$, $\RZ^\uparrow = \D$ and $\Jone^\uparrow =
\J$, and that $\BR_2^\uparrow$ and $\BL_2^\uparrow$ are equal,
respectively, to the pseudovarieties $\RR$ and $\LL$, of $\gR$-trivial
and $\gL$-trivial semigroups.

Then Theorem~\ref{thm: TW} suffices to compute defining
pseudo-identities for all the $\Y^\uparrow$ ($\Y\in\calL(\B)$) -- and
hence to prove the decidability of each of these pseudovarieties.

\begin{example}\label{example G2I2}
    Theorem~\ref{thm: TW} shows that $(\BR'_2 \vee \BL'_2)^\uparrow =
    L\J \cap \DA$ since $ \BR'_2 \vee \BL'_2 = \BR'_3 \cap \BL'_3 =  L(\BR_2 
    \cap \BL_2) \cap \B =  L\Jone \cap \B$ and $\Jone^\uparrow  = \J$.
    
    It also shows that ${\BR'_2}^\uparrow = \RR\cap L\J$, since
    $\BR'_2 = \BR_2 \cap (\BR'_2 \vee \BL'_2)$.
    
    For $\RR = \BR_2^\uparrow$, Theorem~\ref{thm: TW} yields the
    pseudo-identity $x_2^\omega(x_1^\omega x_2^\omega
    x_1^\omega)^\omega = x_2^\omega(x_1^\omega x_2^\omega
    x_1^\omega)^\omega x_2^\omega$.  One can verify that, together
    with the pseudo-identity defining $\DA$, this is equivalent to the
    usual pseudo-identity describing $\RR$, namely $(xy)^\omega =
    (xy)^\omega x$.
    
    The pseudo-identities for the $\BR_m^\uparrow$, $m\ge 3$, are 
    naturally more complicated. For instance, for $m = 3$, we get
    \begin{align*}
	\phi(G_3) &= (x_3^\omega ((x_1^\omega x_2^\omega
	x_1^\omega)^\omega x_2^\omega (x_1^\omega x_2^\omega
	x_1^\omega)^\omega)^\omega x_3^\omega)^\omega \\
	\phi(I_3) &= (x_3^\omega ((x_1^\omega x_2^\omega x_1^\omega)^\omega
	x_2^\omega (x_1^\omega x_2^\omega x_1^\omega)^\omega)^\omega
	x_3^\omega)^\omega \\
	& \qquad \qquad (x_3^\omega ((x_1^\omega x_2^\omega
	x_1^\omega)^\omega x_2^\omega (x_1^\omega x_2^\omega
	x_1^\omega)^\omega)^\omega x_3^\omega)^\omega x_2^\omega
	(x_1^\omega x_2^\omega x_1^\omega)^\omega x_2^\omega.
    \end{align*}
\end{example}

\section{Main result}\label{sec: main result}

Let $m\ge 1$.  It is not difficult to deduce from
Theorems~\ref{lattice bands} and~\ref{thm: TW} that
$\K\malcev\BL_{m}^\uparrow \subseteq \BR_{m+1}^\uparrow$ (and it is
done explicitly in the proof of Theorem~\ref{main theorem} below).  We
prove that the equality actually holds, showing that one can climb in
the lattice $\calL(\DA)$ in a way that directly mimics the steps in
the countable lattice $\calL(\B)$.

\begin{theorem}\label{main theorem}
    For each $m \ge 1$, $\BR_{m+1}^\uparrow = \K \malcev
    \BL_m^\uparrow$ and $\BL_{m+1}^\uparrow = \D \malcev
    \BR_m^\uparrow$
\end{theorem}

\proof
If $m = 1$, the announced equalities are classical results, namely 
the facts that $\RR = \K \malcev \J$ and $\LL = \D \malcev \J$ 
\cite{Pin1986book}. Let us now assume that $m\ge 2$.

If $M$ is a band in $\K\malcev\BL_{m}^\uparrow$, then
Proposition~\ref{congruence malcev} shows that $M/\!\sim_\K \in
\BL_m^\uparrow$.  Since $M/\!\sim_\K$ is a band as well, we have
$M/\!\sim_\K \in \BL_m$.  Moreover, each $\sim_\K$-class is a band,
and a semigroup in $\K$.  Therefore the projection $M \to M/\!\sim_\K$
is an $\LZ$-morphism, and $M \in \LZ\malcev\BL_m = \BR_{m+1}$.  Thus
$(\K\malcev\BL_{m}^\uparrow) \cap \B \subseteq \BR_{m+1}$, and hence
$\K\malcev\BL_{m}^\uparrow \subseteq \BR_{m+1}^\uparrow$.

Conversely, let us assume that $M \in \BR_{m+1}^\uparrow$.  By
Theorem~\ref{thm: TW} (and with the notation in that statement), $M$
satisfies the pseudo-identity $\phi(G_{m+1}) = \phi(I_{m+1})$.  We
want to show that $M/\!\sim_\K \in \BL_m^\uparrow$, that is,
$M/\!\sim_\K$ satisfies $\phi(\overline{G_{m}}) =
\phi(\overline{I_{m}})$.

It is easily verified by induction that the variables which occur in
$G_m$ are the same that occur in $I_m$, namely $x_1,\ldots,x_m$.  We
need to verify that, for each morphism $\psi\colon\{x_1,\ldots,x_m\}^*
\to M$, we have $\psi(\phi(\overline{G_{m}})) \sim_\K
\psi(\phi(\overline{I_{m}}))$.

Let $e\in M$ be an idempotent such that $e\ \psi(\phi(\overline{G_{m}}))
\gJ e$.  Then each $\psi(x_i)$ ($1\le i \le m$) is in $M_e$, the 
subsemigroup of $M$ generated by the elements that are $\gJ$-greater 
than or equal to $e$.

Let us extend $\psi$ to $\{x_1,\ldots,x_{m+1}\}^*$ by letting
$\psi(x_{m+1}) = e$.  Since $\phi(x_{m+1}) = (x_{m+1}^\omega
\phi(\overline{G_{m}}G_{m})^\omega x_{m+1}^\omega)^\omega$ and $eM_ee
= e$ (Proposition~\ref{prop basics DA}), we find that
$\psi(\phi(x_{m+1})) = e$. It follows:
\begin{align*}
    e\ \psi(\phi(\overline{G_{m}})) &=
    \psi(\phi(x_{m+1}\overline{G_{m}})) \\
    &= \psi(\phi(G_{m+1})) \\
    &= \psi(\phi(I_{m+1}))\textrm{ since $M$ satisfies $\phi(G_{m+1})
    = \phi(I_{m+1})$,} \\
    &= \psi(\phi(x_{m+1}\overline{G_m}x_{m+1}\overline{I_m})) \\
    & = e\ \psi(\phi(\overline{G_m}))\ e\
    \psi(\phi(\overline{I_m}))\textrm{ since $\psi(\phi(x_{m+1})) =
    e$,} \\
    & = e\ \psi(\phi(\overline{I_m}))\textrm{ since
    $eM_ee = e$.}
\end{align*}
By symmetry, this shows that $\psi(\phi(\overline{G_{m}})) \sim_\K
\psi(\phi(\overline{I_{m}}))$, which concludes the proof.
\eop

Consequences of this result are explored in the next section.

\section{Applications}\label{sec: applications}

\subsection{Semigroup-theoretic consequences}

An immediate consequence of Theorem~\ref{main theorem} we want to
point out is that we now have explicit pseudo-identities for a number
of natural pseudovarieties.  For instance, no pseudo-identity was
known in the literature for $\K \malcev \LL = \BR_3^\uparrow$ (even
though \cite{PinWeil1996ja} gives general tools to compute this type
of pseudo-identities).  We get
$$\K\malcev\LL = \llbracket \phi(G_3) = \phi(I_3),\ (xy)^\omega x
(xy)^\omega = (xy)^\omega\rrbracket,$$
where $\phi(G_3)$ and $\phi(I_3)$ were computed in
Example~\ref{example G2I2}.

For convenience, in the rest of this paper, we write $\RR_m$ and
$\LL_m$ for $\BR_m^\uparrow$ and $\BL_m^\uparrow$ respectively.  As
indicated in Section~\ref{sec: projecting lattices}, $\RR_1 = \LL_1$
is the pseudovariety of $\gJ$-trivial semigroups, $\RR_2$ is the
pseudovariety of $\gR$-trivial semigroups and $\LL_2$ is the
pseudovariety of $\gL$-trivial semigroups.

We note the following elementary remark (where the case $m = 2$ 
is well-known, see \cite{Pin1986book}).

\begin{proposition}
    For each $m\ge 2$, $\RR_m \cap \LL_m = \Nil\malcev(\RR_m \cap 
    \LL_m)$.
\end{proposition}

\proof
Observe that the operation $\V\mapsto \K\malcev\V$ is idempotent.  In
particular, $\K \malcev (\RR_m \cap \LL_m) \subseteq \K\malcev\RR_m =
\RR_m$.  Similarly, $\D \malcev (\RR_m \cap \LL_m) \subseteq 
\LL_m$, and hence
$$\big(\K \malcev (\RR_m \cap \LL_m)\big) \cap \big(\D \malcev (\RR_m
\cap \LL_m)\big) \subseteq \RR_m \cap \LL_m.$$
The result follows because $(\K\malcev\V) \cap (\D\malcev\V) =
(\K\cap\D)\malcev\V$ \cite[Corollary 3.2]{PinWeil1996ja}, and
$\K\cap\D = \Nil$.
\eop

We now consider the sequences of pseudovarieties $(\RR_m)_m$ and
$(\LL_m)_m$.  It is clear from Theorem~\ref{main theorem} that $\RR_m
\subseteq \RR_{m+2}$, but we have a stronger result.

\begin{proposition}
    For each $m\ge 1$, we have $\RR_m \subseteq \RR_{m+1}$ and $\LL_m
    \subseteq \LL_{m+1}$.  Moreover, both these hierarchies are
    infinite.
\end{proposition}

\proof
By Theorem~\ref{lattice bands}, we have
$$\RR_m\cap\B = \BR_m \subseteq \BR_m\vee\BL_m = 
\BR_{m+1}\cap\BL_{m+1} \subseteq \BR_{m+1}$$
and hence $\RR_m \subseteq \BR_{m+1}^\uparrow = \RR_{m+1}$.  The dual
inclusion, namely $\LL_m \subseteq \LL_{m+1}$ is proved in the same
way.

The infinity of either hierarchy is verified by considering the
sequences $(\RR_m\cap\B)_m = (\BR_m)_m$ and $(\LL_m\cap\B)_m =
(\BL_m)_m$: both these hierarchies are known to be infinite.
\eop

\begin{remark}
    With the same reasoning, one can show that, if $\V$ is a
    pseudovariety containing $\Jone$ and not containing all of $\B$,
    then the sequence of pseudovarieties starting at $\V$ and obtained
    by applying alternately the operations $\X\mapsto\K\malcev\X$ and
    $\X\mapsto\D\malcev\X$ are infinite.
\end{remark}

\begin{proposition}\label{prop: union DA}
    We have $\bigcup_m\RR_m = \bigcup\LL_m = \DA$.  In particular,
    $\DA$ is the least pseudovariety of semigroups containing $\Jone$
    and closed under the operations $\X\mapsto\K\malcev\X$ and
    $\X\mapsto\D\malcev\X$.  In addition, if $M \in \DA$ is
    $m$-generated ($m\ge 2$), then $M$ is in the pseudovariety
    obtained from $\Jone$ by $m$ alternated applications of these
    operations, starting with a Mal'cev product with $\K$ (resp.
    $\D$).
\end{proposition}

\proof
It is immediate that each $\RR_m$ and each $\LL_m$ is contained in
$\DA$.  Conversely, let $M\in \DA$ and let $m\ge 1$ be such that $M$
is $m$-generated.  By Theorem~\ref{thm: TW}~(6), $M \in \RR_{m+1} \cap
\LL_{m+1}$.  This concludes the proof.
\eop

\begin{remark}
    Th\'erien's and Wilke's work \cite{TherienWilke1998stoc} 
    implicitly contains a version of the part of the statement 
    concerning $m$-generated elements of $\DA$, as their proof of 
    the equivalence between $\DA$-recognizability and the 2-variable 
    fragment of first-order logic relies on an induction on the 
    cardinality of the alphabet.
\end{remark}
    
\begin{remark}\label{remark iteration monoids}
    It is well known that $\DA = L\I\malcev\Jone$
    \cite{Schutzenberger1976sf,Pin1986book,Almeida1994book}, and that
    $L\I = \K\vee\D$.  So Proposition~\ref{prop: union DA} states the
    natural-sounding fact that the closure of $\Jone$ under the
    repeated application of the (idempotent) operations
    $\X\mapsto\K\malcev\X$ and $\X\mapsto\D\malcev\X$ is the same as
    its closure under the (idempotent as well) operation $\X\mapsto
    L\I\malcev\X$.  Yet our proof is very specific for $\DA$.  It is
    an interesting question whether this result is in fact more
    general.  We suspect that if $\V\subseteq \DS$ and $\V =
    \K\malcev\V = \D\malcev\V$, then $\V = L\I\malcev\V$, but that
    $\DS$ is the maximal pseudovariety in which this holds.
\end{remark}

\subsection{Language-theoretic consequences}\label{sec: languages}

The language-theoretic corollary we want to record is a simple 
translation of Proposition~\ref{prop: union DA}, but one worth noting.

Recall that if $A$ is an \textit{alphabet} (a finite, non-empty set),
we denote by $A^*$ the free monoid over $A$.  A \textit{language}
$L\subseteq A^*$ is \textit{recognized} by a monoid $M$ if there
exists a morphism $\phi\colon A^*\to M$ such that $L =
\phi\inv(\phi(L))$.  A \textit{class of languages} $\LangV$ is a
collection $\LangV = (\LangV(A))_A$, indexed by all finite alphabets
$A$, such that $\LangV(A)$ is a set of languages in $A^*$.  If $\V$ is
a pseudovariety of monoids, we let $\LangV(A)$ be the set of all
languages of $A^*$ which are recognized by a monoid in $\V$.  The
class $\LangV$ has important closure properties: each $\LangV(A)$ is
closed under Boolean operations and under taking residuals (if $L\in
\LangV(A)$ and $u\in A^*$, then $Lu\inv$ and $u\inv L$ are in
$\LangV(A)$); and if $\phi\colon A^* \rightarrow B^*$ is a morphism
and $L \in \LangV(B)$, then $\phi\inv(L) \in \LangV(A)$.  Classes of
recognizable languages with these properties are called
\textit{varieties} of languages, and Eilenberg's theorem (see
\cite{Pin1986book}) states that the correspondence $\V \mapsto \LangV$,
from pseudovarieties of monoids to varieties of recognizable
languages, is a lattice isomorphism.  Moreover, the decidability of
membership in the pseudovariety $\V$, implies the decidability of the
variety $\LangV$: indeed, a language is in $\LangV$ if and only if its
(effectively computable) syntactic monoid is in $\V$.

Let $K,L$ be languages in $A^*$ and let $a\in A$.  The product $KaL$
is said to be \textit{deterministic} if each word $u\in KaL$ has a
unique prefix in $Ka$.  If $k \ge 1$, $L_0,\ldots,L_k$ are languages
in $A^*$ and $a_1,\ldots,a_k \in A$, the product $L_0a_1L_1\cdots
a_kL_k$ is said to be \textit{deterministic} if the products $L_{i-1}
a_i (L_i\cdots a_nL_n)$ are deterministic, for $1 \le i \le n$.

Dually, the product $KaL$ is said to be \textit{co-deterministic} if
each word $u\in KaL$ has a unique suffix in $aL$.  The product
$L_0a_1L_1\cdots a_kL_k$ is said to be \textit{co-deterministic} if
the products $(L_0 a_1 \cdots L_{i-1})a_{i}L_i$ are co-deterministic,
for $1 \le i \le n$.
 
Finally, the product $L_0a_1L_1\cdots a_kL_k$ is said to be
\textit{unambiguous} if every word $u$ in this language admits a
unique decomposition in the form $u = u_0a_1u_1 \cdots a_nu_n$ with
each $u_i \in L_i$.  It is easily verified that a deterministic or
co-deterministic product is a particular case of an unambiguous
product.

These operations are extended to classes of languages: If $\LangV$ is
a class of languages, let $\LangV^{det}$ (resp.  $\LangV^{codet}$)
denote the class of languages such that, for each alphabet $A$,
$\LangV^{det}(A)$ (resp.  $\LangV^{codet}(A)$) is the set of all
Boolean combinations of languages of $\LangV(A)$ and of deterministic
(resp.  co-deterministic) products of languages of $\LangV(A)$.  Let
also $\LangV^{unamb}$ be the class of languages such that, for each
alphabet $A$, $\LangV^{unamb}(A)$ is the set of all finite unions of
unambiguous products of languages of $\LangV(A)$.

Sch\"utzenberger \cite{Schutzenberger1976sf,Pin1986book} gave algebraic characterizations
of the closure operations $\LangV \longmapsto \LangV^{det}$, $\LangV
\longmapsto \LangV^{codet}$ and $\LangV \longmapsto \LangV^{unamb}$
for varieties of languages: he showed that $\LangV^{det}$, 
$\LangV^{codet}$ and $\LangV^{unamb}$ are varieties of languages, and 
the the corresponding pseudovarieties are $\K\malcev\V$, 
$\D\malcev\V$ and $L\I\malcev\V$, respectively.

Proposition~\ref{prop: union DA} now easily translates to the
following statement.

\begin{proposition}\label{prop: closure products}
    The least variety of languages containing the languages of the
    form $B^*$ ($B \subseteq A$) and closed under deterministic and
    co-deterministic product, is the variety corresponding to $\DA_\M$.
    
    More precisely, every unambiguous product of languages
    $B_1^*a_1B^*_2\cdots a_kB^*_{k+1}$ where the $B_i$ are subsets of
    alphabet $A$, can be expressed in terms of the $B^*_i$ and the
    $a_i$ using only Boolean operations and at most $|A|$ alternated
    applications of the deterministic and co-deterministic products 
    -- starting with a deterministic (resp. co-deterministic) product.
\end{proposition}

\begin{remark}\label{remark iteration languages}
    As in Remark~\ref{remark iteration monoids}, it is interesting to
    note that, while this result (that unambiguous products can be
    expressed by iterated deterministic and co-deterministic products)
    sounds natural, its proof is very specific for $\Jone$ and $\DA$:
    does it hold in general?
\end{remark}

%
%
%
%
	
\bibliographystyle{plain}
{\small


\newcommand{\Ju}{Ju}\newcommand{\Th}{Th}\newcommand{\Yu}{Yu}

}

\end{document}